\documentclass[11pt]{article}
\usepackage[numbers,sort&compress]{natbib}
\usepackage[colorlinks,
            linkcolor=blue,
            anchorcolor=green,
            citecolor=magenta
            ]{hyperref}
\usepackage{mathrsfs,amssymb,amsfonts}
\usepackage{graphicx}
\usepackage{amsmath,amssymb,latexsym,color}
\usepackage[mathscr]{eucal}
\usepackage[numbers]{natbib}
\newtheorem{thm}{Theorem}[section]

\newtheorem{prop}[thm]{Proposition}
\newtheorem{lem}[thm]{Lemma}

\newtheorem{example}{Example}
\newtheorem{remark}{Remark}[section]

\newcommand{\rc}{{\rm rc}}

\newcommand{\proof}{{\it Proof.\quad}}
\newcommand{\qed}{\hfill\Box\medskip}

\usepackage{CJK}
\begin{document}
\begin{CJK*}{GBK}{song}
\renewcommand{\abovewithdelims}[2]{
\genfrac{[}{]}{0pt}{}{#1}{#2}}

\title{\bf Rainbow connectivity of the non-commuting graph of a finite group}

\author{Yulong Wei\quad Xuanlong Ma\quad Kaishun Wang\\
{\footnotesize   \em  Sch. Math. Sci. {\rm \&} Lab. Math. Com. Sys., Beijing Normal University, Beijing, 100875,  China} }
 \date{}
 \date{}
 \maketitle

\begin{abstract}
Let $G$ be a finite non-abelian group.
The non-commuting graph $\Gamma_G$ of $G$ has the vertex set $G\setminus Z(G)$ and two distinct vertices $x$ and $y$ are adjacent if
$xy\ne yx$, where $Z(G)$ is the center of $G$.
We prove that the rainbow $2$-connectivity of $\Gamma_G$ is $2$. In particular, the rainbow connection number of $\Gamma_G$ is $2$.
Moreover, for any positive integer $k$, we prove that there exist infinitely many non-abelian groups $G$ such that the rainbow $k$-connectivity of $\Gamma_G$ is $2$.


\footnote{E-mail addresses: weiyl@mail.bnu.edu.cn (Y. Wei), xuanlma@mail.bnu.edu.cn (X. Ma),\\ wangks@bnu.edu.cn (K. Wang).}
\end{abstract}

\section{Introduction}

Let $\Gamma$ be a connected graph with the vertex set $V(\Gamma)$ and the  edge set $E(\Gamma)$.
Given an edge coloring of $\Gamma$.
A path $P$ is {\em rainbow} if no two edges of $P$ are colored the same.
The {\em vertex connectivity} of $\Gamma$, denoted by $\kappa(\Gamma)$, is the smallest number of vertices whose deletion from $\Gamma$ disconnects it.
For any positive integer $k\le \kappa(\Gamma)$,
an edge-colored graph is called
{\em rainbow-$k$-connected} if any two distinct vertices of $\Gamma$ are
connected by at least $k$ internally disjoint rainbow
paths.
The {\em rainbow-$k$-connectivity} of $\Gamma$, denoted by $\rc_k(\Gamma)$,
is the minimum number of colors required to color the
edges of $\Gamma$ to make it rainbow-$k$-connected.
We usually denote $\rc_1(\Gamma)$ by rc$(\Gamma)$, which is called the {\em rainbow connection number} of $\Gamma$.

In \cite{CJMZ} and \cite{CJMZh}, Chartrand et al. first introduced the concept of rainbow $k$-connectivity for $k=1$ and $k\ge 2$, respectively.
Rainbow k-connectivity has application in transferring information of high security in communication networks. For details we refer to \cite{CJMZh} and \cite{Eri}.
The NP-hardness of determining $\mathrm{rc(\Gamma)}$ was shown by Chakraborty et al. \cite{CFMY}. Recently, the rainbow connectivity of some special classes of graphs have been studied;
see \cite{LSue12} for complete graphs, \cite{LSue11} for regular complete bipartite graphs,
\cite{LLL,LMa1,LMa} for Cayley graphs and \cite{MFW} for power graphs.
For more information, see \cite{LSus}.

For a non-abelian group $G$, the {\it non-commuting graph} $\Gamma_G$ of $G$ has the vertex set $G\setminus Z(G)$ and two distinct vertices $x$ and $y$ are adjacent if
$xy\ne yx$, where $Z(G)$ is the center of $G$.
According to \cite{Ne} non-commuting graphs were first considered by Erd\H{o}s in 1975.
Over the past decade, non-commuting graphs have received considerable attention. For example,
Abdollahi et al. \cite{AAM} proved that the diameter of any non-commuting graph is $2$.
For two non-abelian groups with isomorphic non-commuting graphs,
the sufficient conditions that guarantee their orders are equal were provided
by Abdollahi and Shahverdi \cite{AS} and Darafsheh \cite{Da}. Akbari and  Moghaddamfar \cite{AM} studied strongly regular non-commuting graphs.  Solomon and Woldar \cite{SW} characterized some simple groups
by their non-commuting graphs.

In this paper we study the rainbow $k$-connectivity of non-commuting graphs and obtain
the following results.

\begin{thm}\label{main2}
Let $G$ be a finite non-abelian group. Then $\rc_2(\Gamma_G)=2$. In particular, $\rc(\Gamma_G)=2$.
\end{thm}

\begin{thm}\label{main3}
For any positive integer $k$, there exist infinitely many non-abelian groups $G$ such that
$\rc_k(\Gamma_G)=2$.
\end{thm}

\section{Preliminaries}
In this section we present some lemmas which we need in the sequel.

For vertices $x,y$ of a graph $\Gamma$, let
$\tau(x,y)$ be the number of the common neighbors of $x$ and $y$.

\begin{lem}\label{db} Let $G$ be a finite non-abelian group, and let $x$ and $y$ be two  distinct vertices of $\Gamma_G$. Then
$\tau(x,y)\ge \frac{1}{6}|G|$.
\end{lem}
\proof For each $g\in G$, $C_G(g)$ denotes the centralizer of $g$ in $G$.
By the principle of inclusion and exclusion,
\begin{equation*}\label{cap}
\begin{aligned}
\tau(x,y)&=|G|-|C_G(x)\cup C_G(y)|\\
&=|G|-|C_G(x)|-|C_G(y)|+|C_G(x)\cap C_G(y)|\\
&\ge |G|-|C_G(x)|-|C_G(y)|+\frac{|C_G(x)|\cdot|C_G(y)|}{|G|}.
\end{aligned}
\end{equation*}
If $|C_G(x)|=|C_G(y)|=\frac{1}{2}|G|$, then
$$\tau(x,y)\geq \frac{|C_G(x)|\cdot|C_G(y)|}{|G|}=\frac{1}{4}|G|;$$
if not, $$~~~~~~~~~~~~~~~~~~~~~~~~~~~~\tau(x,y)\geq |G|-|C_G(x)|-|C_G(y)|\geq \frac{1}{6}|G|.~~~~~~~~~~~~~~~~~~~~~~\qed $$

The lexicographic product $\Gamma \circ \Lambda$ of graphs $\Gamma$ and $\Lambda$ has the vertex set $V(\Gamma)\times V(\Lambda)$, and two vertices $(\gamma, \lambda), (\gamma' , \lambda')$ are adjacent if $\{\gamma,\gamma'\}\in E(\Gamma)$, or if $\gamma=\gamma'$ and $\{\lambda,\lambda'\}\in E(\Lambda)$.

\begin{lem}\label{zdj}
Let $G$ be a non-abelian group and $A$ be an abelian group of order $n$. Then
$$\Gamma_{G\times A}\cong \Gamma_G \circ \overline{K_{n}},$$
where $\overline{K_{n}}$ is the complement of the complete graph $K_n$.
\end{lem}

For positive integers $l,r$ and $t$,
let $K_{l[r]}$ denote a complete $l$-partite graph with each part of order $r$, and let
$K_{l[r],t}$ denote a complete $(l+1)$-partite graph with $l$ parts of order $r$ and a part of order $t$.

\begin{lem}\label{dg} Let $D_{2n}$ and $Q_{4m}$ be respectively the dihedral group of order $2n$ and the generalized quaternion group of order $4m$, where $n\ge 3$ and $m\ge 2$. Then

${\rm (i)}$ If $n$ is odd, then $\Gamma_{D_{2n}}\cong K_{n[1],n-1}$.

${\rm (ii)}$ If $n$ is even, then $\Gamma_{D_{2n}}\cong  K_{\frac{n}{2}[2],n-2}$.

${\rm (iii)}$ $\Gamma_{Q_{4m}}\cong  K_{m[2],2m-2}$.
\end{lem}

Li and Sun \cite{LSue12} studied the rainbow $k$-connectivity of some families of  complete multipartite graphs.
Now we compute the rainbow $k$-connectivity of another family.

\begin{prop}\label{mainlem}
Let $m\ge n+1$, $lmn\neq2$. Then
$\rc_2(K_{m[l],ln})=2$.
\end{prop}
\proof
Write $\Gamma=K_{m[l],ln}$. Let $\{a_{j,i}: 1\le j \le l\}$ and $\{a_{j,m+1}: 1\le j \le ln\}$ be all parts of $\Gamma$, where $i=1,\ldots,m$.

\medskip

{\em Case 1.} $n=1$.

{\em Case 1.1.}  $m=2$.

If $l=2r$, then we assign a color to the edges
\begin{equation}\label{ca1}
\begin{array}{l}
\{a_{2j-1,1},a_{2j-1,2}\},\{a_{2j-1,1},a_{2j,2}\},\{a_{2j-1,1},a_{2j-1,3}\},
\{a_{2j,1},a_{2j-1,2}\}, \\
\{a_{2j,1},a_{2j,2}\}, \{a_{2j,1},a_{2j,3}\},\{a_{2j-1,2},a_{2j-1,3}\},\{a_{2j-1,2},a_{2j,3}\}, ~~1\le j \le r
\end{array}
\end{equation}
and another color to the remaining edges.

If $l=2r+1$, then we assign a color to the edges
\begin{equation*}\label{ca2}
\begin{array}{l}
\{a_{l,2},a_{l,3}\}, \{a_{l,1},a_{2j-1,2}\},\{a_{l,1},a_{2j,3}\},\{a_{l,2},a_{2j-1,1}\},
\{a_{l,2},a_{2j,1}\},\\
\{a_{l,2},a_{2j-1,3}\}, \{a_{l,2},a_{2j,3}\},\{a_{l,3},a_{2j,2}\},~~1\le j \le r
\end{array}
\end{equation*}
and the edges in (\ref{ca1}), and another color to all other edges.

{\em Case 1.2.}  $m=3$.

The edges $$\{a_{j,1},a_{j,2}\},\{a_{j,2},a_{j,4}\},\{a_{j,3},a_{j,4}\},~~1\le j \le l$$ are assigned by a color and all other edges are assigned by another color.

{\em Case 1.3.}  $m\ge 4$.

The edges
$$\{a_{j,i},a_{j,i+1}\},\{a_{j,m+1},a_{j,1}\},~~1\le i \le m,~~1\le j \le l$$
are assigned by a color and all other edges are assigned by another color.

Note that all the above colorings make $\Gamma$ rainbow-$2$-connected. Hence $\rc_2(\Gamma)=2$.
\medskip

{\em Case 2.} $n\ge 2$.

We assign a color to
$$\{a_{j,i},a_{j,i+1}\},\{a_{j,1},a_{j,m}\},\{a_{j,i'},a_{(j-1)n+i',m+1}\}, 1\le i\le m-1,1\le i' \le n,1\le j \le l$$
and  another color to the remaining edges.
Note that this coloring makes $\Gamma$ rainbow-$2$-connected.
This implies that $\rc_2(\Gamma)=2$.
$\qed$

\section{Proof of main results}\label{rcd}
In this section, we shall prove  Theorems \ref{main2} and \ref{main3}.

\begin{prop}\label{p1}
Let $G$ be a finite non-abelian group with $|G|\ge 114$. Then $\rc_2(\Gamma_G)=2$.
\end{prop}
\proof
We randomly color the edges of $\Gamma_G$ with two colors. Denote by $\mathcal{P}_G$ the probability that such a random coloring makes it not rainbow-$2$-connected.
It suffices to prove that $\mathcal{P}_G<1$.

Let $x$ and $y$ be two distinct vertices of $\Gamma_G$.
If $x$ and $y$ are adjacent, then the probability that there exist no
rainbow paths of length $2$ from $x$ to $y$ is $(1/2)^{\tau(x,y)}$.
If  $x$ and $y$ are non-adjacent, then the probability that $\Gamma_G$ has precisely a rainbow path of length $2$ from $x$ to $y$ is $\tau(x,y)(1/2)^{\tau(x,y)}$, and the  probability that $\Gamma_G$ has no rainbow paths of length $2$ from $x$ to $y$ is $(1/2)^{\tau(x,y)}$.
Note that
\begin{equation}\label{en}
|E(\Gamma_G)|=\frac{1}{2}\sum_{x\in V(\Gamma_G)}(|G|-|C_G(x)|)
\ge \frac{1}{4}|G|(|G|-|Z(G)|).
\end{equation}
Write
\begin{equation}\label{pp}
\mathcal{P}=\sum_{x\sim y}\left(\frac{1}{2}\right)^{\tau(x,y)}+
\sum_{x\nsim y}\left(\left(\frac{1}{2}\right)^{\tau(x,y)}+\tau(x,y)\left(\frac{1}{2}\right)^{\tau(x,y)}\right),
\end{equation}
where $x\sim y$ denotes that $x$ and $y$ are adjacent.
Now
\begin{eqnarray*}
\mathcal{P}_G&\le & \mathcal{P}\\
&\le& \sum_{x\sim y}\left(\frac{1}{2}\right)^{\frac{1}{6}|G|}
+\sum_{x\nsim y}\left(\frac{1}{2}\right)^{\frac{1}{6}|G|}+\sum_{x\nsim y}|G|\left(\frac{1}{2}\right)^{\frac{1}{6}|G|}  ~~~~~~~~({\rm by~Lemma~\ref{db}}) \\
&=& \binom{|V(\Gamma_G)|}{2}\left(\frac{1}{2}\right)^{\frac{1}{6}|G|}+\sum_{x\nsim y}|G|\left(\frac{1}{2}\right)^{\frac{1}{6}|G|}\\
&\le & \binom{|G|-|Z(G)|}{2}\left(\frac{1}{2}\right)^{\frac{1}{6}|G|}\\
&&+|G|\left(\frac{1}{2}\right)^{\frac{1}{6}|G|}\left(\binom{|G|-|Z(G)|}{2}-
\frac{1}{4}(|G|-|Z(G)|)|G|\right)   ~~~~~~~~({\rm by}~ (\ref{en})) \\
&=&\frac{1}{4}(|G|-|Z(G)|)(|G|^2-2|Z(G)|-|Z(G)||G|-2)\left(\frac{1}{2}\right)^{\frac{1}{6}|G|}\\
&<&\left(\frac{1}{2}\right)^{\frac{1}{6}|G|+2}|G|^3.
\end{eqnarray*}
It follows that if $|G|\ge 114$, then $\mathcal{P}_G<1$. $\qed$

\begin{prop}\label{p2}
Let $G$ be a finite non-abelian group with $|G|< 114$. Then $\rc_2(\Gamma_G)=2$.
\end{prop}
\proof
Let $\mathcal{P}_G$ be the probability that such a random coloring makes $\Gamma_G$ not rainbow-$2$-connected. Thus if $\mathcal{P}_G\le 1$, then $\rc_2(\Gamma_G)=2$.
Using GAP \cite{G13}, we compute $\mathcal{P}$ (see (\ref{pp}))
by the following code.
\begin{verbatim}
M:=Elements(G);
k:=0;
s:=0;
for i in [1..Size(M)] do
 if Centralizer(G,M[i])<>G then
  for j in [1..Size(M)] do
   if Centralizer(G,M[j])<>G and IsAbelian(Group(M[i],M[j]))=false
                                                and M[i]<>M[j] then
               t:=Order(G)-Size(Union(Elements(Centralizer(G,M[i])),
                                    Elements(Centralizer(G,M[j]))));
               k:=k+(1/2)^t;
   fi;
  od;
 fi;
od;
for i in [1..Size(M)] do
 if Centralizer(G,M[i])<>G then
  for j in [1..Size(M)] do
   if Centralizer(G,M[j])<>G and IsAbelian(Group(M[i],M[j]))
                                         and M[i]<>M[j] then
       t:=Order(G)-Size(Union(Elements(Centralizer(G,M[i])),
                            Elements(Centralizer(G,M[j]))));
       s:=s+(1/2)^t+t*(1/2)^t;
   fi;
  od;
 fi;
od;
P:=(s+k)/2;
\end{verbatim}
By this code, one gets that $\mathcal{P}_G\le \mathcal{P}<1$ with the following exceptions:

(i) $D_6, D_8, Q_8, D_{10}, D_{12}, Q_{12}, D_{14}, D_6\times \mathbb{Z}_3,
D_8\times \mathbb{Z}_3, Q_8\times \mathbb{Z}_3$, and all non-abelian groups of order $16$.

(ii) $G_1$ and $G_2$, where $G_1=$ SmallGroup(32,49) and $G_2=$ SmallGroup(32,50) in GAP.

Note that $\Gamma_H\cong K_{4[2],6}$ or $K_{3[4]}$ for any non-abelian group $H$ of order $16$. By Lemmas \ref{zdj}, \ref{dg} and Proposition \ref{mainlem}, the rainbow $2$-connectivity of the non-commuting graph of each group in (i) is $2$.

Next we shall prove that $\rc_2(\Gamma_{G_1})=\rc_2(\Gamma_{G_2})=2$.
Note that $\Gamma_{G_1}\cong J(6,2)\circ \overline{K_2}$.
For convenience, in the following we use $ab$ to denote the set $\{a,b\}$ for  two distinct letters $a,b$. Write $\Gamma=J(6,2)\circ \overline{K_2}$ and
$$V(\Gamma)=\{a_ia_j,b_ib_j:1\le i,j\le 6,i\ne j\},$$
$$E(\Gamma)=\{\{a_ia_j,a_ia_k\},\{b_ib_j,b_ib_k\},\{a_ia_j,b_ib_k\}:1\le i,j,k\le 6,i\ne j,i\ne k, j\ne k\}.$$
We assign the edges in $\{\{a_ia_j,a_ia_k\}: i> \max\{j,k\}\}$, $\{\{b_ib_j,b_ib_k\}: i> \max\{j,k\}\}$ and $\{\{a_ia_j,b_ib_k\}: i< \min\{j,k\}\}$  a color and all other edges another color.
It follows that $\rc_2(\Gamma)=2$. Hence $\rc_2(\Gamma_{G_1})=\rc_2(\Gamma_{G_2})=2$.
$\qed$

Combining Propositions \ref{p1} and \ref{p2}, we complete the proof of Theorem \ref{main2}.

\bigskip

\noindent{\em Proof of Theorem \ref{main3}:}
By Theorem \ref{main2}, we may assume $k\ge 3$.
Note that for any non-abelian group $H$, $\kappa(\Gamma_H)$ is divisible by $|Z(H)|$ by \cite[Proposition 2.4]{AAM}. Therefore, we may
choose a non-abelian group $G$ with $k\le \kappa(\Gamma_G)$.
Next we prove that $\rc_k(\Gamma_G)=2$  if  $|G|$ is large enough.

We randomly color the edges of $\Gamma_G$ with two colors.
Denote by $\mathcal{P}_G$
the probability that such a random coloring makes it not rainbow-$2$-connected.
It suffices to show that $\mathcal{P}_G<1$.
Let $x$ and $y$ be distinct vertices of $\Gamma_G$.
If $x$ and $y$ are adjacent, then the probability that there exist no
$k$ rainbow paths of length $2$ from $x$ to $y$ is
$$\sum_{i=0}^{k-2}\binom{\tau(x,y)}{i}(1/2)^{\tau(x,y)}.$$
If $x$ is not adjacent to $y$, then the probability that there are no
$k$ rainbow paths of length $2$ from $x$ to $y$ is
$$\sum_{i=0}^{k-1}\binom{\tau(x,y)}{i}(1/2)^{\tau(x,y)}.$$
Write $|G|=n$.
Then we have
\begin{eqnarray*}
\mathcal{P}_G&\le &\sum_{x\sim y}\sum_{i=0}^{k-2}\binom{\tau(x,y)}{i}\left(\frac{1}{2}\right)^{\tau(x,y)}+
\sum_{x\nsim y}\sum_{i=0}^{k-1}\binom{\tau(x,y)}{i}\left(\frac{1}{2}\right)^{\tau(x,y)}\\
&\le&\sum_{x\sim y}\sum_{i=0}^{k-2}\tau(x,y)^i\left(\frac{1}{2}\right)^{\tau(x,y)}+
\sum_{x\nsim y}\sum_{i=0}^{k-1}\tau(x,y)^i\left(\frac{1}{2}\right)^{\tau(x,y)}\\
&\le& \left(\frac{1}{2}\right)^{\frac{1}{6}n}\left(\sum_{x\sim y}\sum_{i=0}^{k-2}n^i+
\sum_{x\nsim y}\sum_{i=0}^{k-1}n^i\right) ~~~~~~~~~~~~~~~~~~~~~~~~~~~~~({\rm by~ Lemma}~\ref{db})\\
&=& \left(\frac{1}{2}\right)^{\frac{1}{6}n}\left(\sum_{i=0}^{k-2}\binom{|V(\Gamma_G)|}{2}n^i+
\left(\binom{|V(\Gamma_G)|}{2}-|E(\Gamma_G)|\right)n^{k-1}\right)\\
&<& \left(\frac{1}{2}\right)^{\frac{1}{6}n}\left(\sum_{i=0}^{k-2}n^{i+2}+n^{k+1}\right)\\
&=&\frac{\sum_{i=2}^{k+1}n^{i}}{2^{\frac{n}{6}}}.
\end{eqnarray*}
This implies that $\mathcal{P}_G<1$ if $n$ is large enough. $\qed$

\section*{Acknowledgement}
This research is supported by National Natural Science Foundation of China (11271047,
11371204) and the Fundamental Research Funds for the Central University of China.

\end{CJK*}

\end{document}